\documentclass[a4paper,12pt]{article}
\usepackage{amsmath,amssymb}
\usepackage{hyperref}

\begin{document}

\title{ON SHIMURA SUBVARIETIES GENERATED BY FAMILIES OF ABELIAN COVERS OF $\mathbb{P}^{1}$}
\author{Abolfazl Mohajer, Kang Zuo}
\date{February 2014}
\maketitle

\

\

\begin{abstract}

We study the locus of Jacobians of  abelian Galois covers of
$\mathbb{P}^{1}$ in $A_{g}$. We find examples of Shimura (special)
subvarieties generated by families of these covers in the Torelli
locus $T_g$ inside $A_g$ and give classifications of such families
giving rise to Shimura subvarieties in $T_g$. Our methods are
based on Moonen-Oort works as well as characteristic $p$
techniques of Dwork and Ogus and Monodromy computations.

\end{abstract}

\

\

\section{Introduction}

\

\

Let $M_{g}$ be the (coarse) moduli space of non-singular complex
curves of genus $g$ and $A_{g}$ the moduli space of principally
polarized complex abelian varieties of dimension $g$. Let
$j:M_{g}\rightarrow A_{g}$ be the Torelli map, i.e., the map which
associates to a curve its Jacobian. The Torelli theorem asserts
that this map is injective. We call the image of this map the
\emph{open Torelli locus} and denote it by $T^{\circ}_{g}$ . By
its definition, it consists of Jacobians of smooth curves of genus
$g$. The closure of the open Torelli locus inside $A_{g}$ is
called the \emph{closed Torelli locus}, which we denote by $T_{g}$
. The closed Torelli locus is the locus of Jacobians of stable
curves of compact type. The quest of Shimura subvarieties (or
special subvarieties) contained generically in the Torelli locus
$T_{g}$ and not fully contained in the boundary is a longstanding
problem, which can be traced back at least to Shimura (see
\cite{63}). Although not realized then, more than two decades
later, Coleman formulated his famous conjecture suggesting that
for $g\geq 4$ there are only finitely many curves $C$  of genus
$g$ over the complex numbers whose Jacobian $J(C)$ is an abelian
variety of CM type. See \cite{87}. Shortly after that, de Jong and
Noot disproved that conjecture by finding examples of families of
curves which give rise to Shimura subvarieties in $A_{g}$ lying
generically in the (open) Torelli locus and intersecting the
(closed) Torelli locus non-trivially. The families that they found
- which were also found in the aforementioned article of Shimura-
were all families of cyclic coverings of $\mathbb{P}^{1}$.  Their
examples were of fiber genus $4$ and $6$ and they did explain the
relation between their examples and the Coleman conjecture.
Examples of higher genera ($g=5$ and $g=7$), again of families of
cyclic coverings of $\mathbb{P}^{1}$, were found by Rohde
\cite{09}. Finally Moonen completed the list of Shimura
subvarieties arising from families of cyclic coverings of
$\mathbb{P}^{1}$ in \cite{11} by showing that there are no more
Shimura subvarieties in this locus. The fiber genus of all of
these families were bounded by $8$. Note that all of these
suvarieties are of PEL type. For some non-PEL type examples in
$A_g$ see for example \cite{17}. Based on the above examples, one
can reformulate the Coleman conjecture as follows:

\

\

\textbf{The (corrected) Coleman conjecture}. For $g\geq 8$ there
are only finitely many smooth projective curves $C$ over
$\mathbb{C}$ of genus $g$ such that $J(C)$ is an abelian variety
of CM type.

\

\

Bearing in mind that every Shimura variety contains a dense subset
of CM points, we see that the Coleman conjecture is related to the
following conjecture and in fact disproving this conjecture will
disprove the Coleman conjecture:

\

\

\textbf{Conjecture (\cite{97}, \S5).} For large $g$ (in any case
$g\geq 8$), there is no positive-dimensional Shimura subvariety
contained in $T_{g}$ such that $Z\cap T^{\circ}_{g}$ is non-empty.

\

\

We just remark that if one assumes the Andr\'e-Oort conjecture to
be true then the above conjecture is even equivalent to the
Coleman conjecture. Note that the condition that our subvariety
meets $T^{\circ}_{g}$ is a key condition. Otherwise one can easily
construct a lot of Shimura subvarieties contained fully in the
boundary of $T_{g}$ for every $g\geq 2$. There are also several
other results on the occurrence of Shimura subvarieties in the
Torelli locus. Viehweg and Zuo for example studied in \cite{04}
the occurrence of Shimura curves in the moduli stack of
principally polarized abelian varieties relating it to Arakelov
equalities. See also \cite{06}. The recent work of Lu and Zuo
\cite{13} has made some progress in the proof of Coleman's
conjecture. They show that there do  not exist Shimura curves in
Torelli locus of hyperelliptic curves of genus greater than $7$.
They also show that Shimura curves with maximal Higgs field (they
are either self product of universal family of elliptic curves or
Shimura curves of Mumford type) can not be contained in Torelli
locus of curves of genus greater than $5$. In \cite{10} Oort and
Moonen asked whether one can obtain further Shimura subvarieties
in the Torelli locus by taking families of abelian coverings of
$\mathbb{P}^{1}$ with a non-cyclic Galois group. They also gave
some examples of families of abelian covers of $\mathbb{P}^{1}$
which give rise to Shimura subvarieties in $T_{g}$.

\

\textbf{Structure of the paper and main results.} In this article
we try to generalize the methods of Moonen-Oort and classify
Shimura subvarieties arising from families of abelian covers of
the projective line. In section $2$ we describe the construction
of abelian covers of $\mathbb{P}^{1}$ and their families. More
precisely, we fix integers $N\geq 2$, $s\geq 4$, $m\geq 1$ and an
$m \times s$ matrix $A$ with entries in
$(\mathbb{Z}/N\mathbb{Z})$. Given an $s$-tuple
$t=(z_{1},...,z_{s})\in (\mathbb{A}^{1}_{\mathbb{C}})^{s}$ there
is a Galois cover $Y_{t}\rightarrow \mathbb{P}^{1}$ branched at
the points $z_{j}$ with local monodromies encoded in the matrix
$A$ and an abelian Galois group $G$ which is isomorphic to the
column span of the matrix $A$ and hence is a subgroup of the group
$(\mathbb{Z}/N\mathbb{Z})^{m}$. By varying the branch points we
obtain a family of abelian covers of $\mathbb{P}^{1}$. The
corresponding family of Jacobians gives a subvariety of $A_{g}$
whose closure $Z=Z(N,s,A)\subset A_{g}$ lies in $T_{g}$. Of
course, $Z$ will be of dimension $s-3$ where $s$ is the number of
branch points of the covering and it lies in the Torelli locus
cutting the open Torelli locus non-trivially. Our goal is to
classify the cases where $Z$ is a Shimura subvariety. We observe
that there are two Shimura subvarieties containing $Z$. Our method
here is a generalization of that of Moonen-Oort \cite{10} and
Moonen \cite{11}. The Jacobians in our families admit an action of
the group ring $\mathbb{Z}[G]$. This action defines a Shimura
subvariety $S(G)$ in $A_{g}$. Note that $S(G)$ does depend also on
the matrix $A$ but since in our construction in section 1 the
matrix, up to possible permutations of rows and columns which
yield isomorphic families, is fixed for each family from the
beginning, we omit this dependence and use the notation $S(G)$.
This subvariety is in fact the Hodge locus given by the Hodge
classes that are exactly elements of $\mathbb{Z}[G]$ viewed as
endomorphisms of the Jacobians in our family. It follows that
$Z\subseteq S(G)$ and therefore $s-3 \leq \dim S(G)$. When $\dim
S(G)=s-3$, or equivalently when $Z=S(G)$, one concludes that the
subvariety $Z$ is a Shimura subvariety. The dimension of $S(G)$
can be calculated in terms of dimensions of certain eigenspaces
with respect to the Galois group action and the equality $\dim
S(G)=s-3$ is then easy to check by a computer and by using a
simple computer program we list the examples in a table (see Table
1). Note that this table appears also in \cite{10}, where they
present these families as a set of examples of families of abelian
covers which give rise to Shimura subvarieties. Their method
-based on the only example that they give- is to decompose the
Jacobian of the fibers by representing them as cyclic covers and
then compute the multiplicity of the corresponding group action on
the tangent space of the Jacobian. Also, they do not claim that
this table contains all examples for which $\dim S(G)=s-3$ and do
not claim that their table is exhaustive (unlike Table $1$ in
\cite{11} which is shown there to be exhaustive). Here we use a
systematic approach based on an explicit formula in terms of the
dimension of the eigenspaces and then we compute also these
dimensions explicitly. Of course, when $s-3< \dim S(G)$ one can
not conclude that $Z$ is not a Shimura subvariety. It can very
well happen that $Z$ is still a smaller Shimura subvariety
contained in $S(G)$. The second Shimura subvariety containing $Z$
is the smallest Shimura subvariety with this property and we
denote it by $S_f$. Its corresponding reductive group $M$ is the
generic Mumford-Tate group of the family. Therefore $Z=S_f$ if and
only if $Z$ is a Shimura subvariety, see Construction 2.2.2 where
we describe how to compute $\dim S_f$. The above equivalence means
that $Z$ is a Shimura subvariety if and only if $\dim S_f=s-3$. We
will use these facts mainly in section 6. In section 3, we state
our first main result, namely, Theorem 3.1. This Theorem asserts
that Table 1 contains various examples of matrices $A$ such that
$Z$ is a Shimura subvariety and moreover: If $s=4$ and $\dim
S(G)=1$, then $A$ shows up in this list. If $s\in \{5,6,7\}$,
$2\leq m \leq 5$ and $N\leq 20$, then $A$ appears in this list. In
\cite{11}, Moonen shows that when $G$ is cyclic, then $Z$ is a
Shimura subvariety if and only if $Z=S(G)$. As described above,
this is by no means obvious. Our objective in the rest of the
paper is to show this equivalence for $G$ abelian. In order to do
this, we require tools from characteristic $p$. Section 4 and 5
provide us with these technical tools. More precisely, in section
4 we describe an obstruction of Dwork-Ogus which vanishes if the
family gives rise to a Shimura subvariety in $A_g$. In section 5
we give a formula for the Hasse-Witt map of the family which will
be used to show that the obstruction does not vanish for our
families. This is done in section 6.1. Our second main theorem,
Theorem 6.1.1, shows that if $s=4$ and $Z$ is a Shimura subvariety
then $Z=S_f=S(G)$ and hence it appears in Table 1. Taken together
with Theorem 3.1, these show that Table 1 completely classifies
one-dimensional Shimura subvarieties which arise from irreducible
abelian covers of $\mathbb{P}^{1}$. Finally, our last main result,
Theorem 6.2.6, shows that when $s$ is large enough ($s>19$), then
the families of abelian covers of $\mathbb{P}^{1}$ do not give
rise to Shimura subvarieties. We achieve this by showing that in
these cases $s-3<\dim S_f$, using constructions form section 2.

\

\

 \section{Construction of abelian covers of $\mathbb{P}^{1}$ and
     their families}
% %\label{}

% % bibitems, please use
%%


\

\

An abelian Galois cover of $\mathbb{P}^{1}$ is determined by a
collection of equations in the following way: Consider an $m\times
s$ matrix $A=(r_{ij})$ whose entries $r_{ij}$ are in
$\mathbb{Z}/N\mathbb{Z}$ for some $N\geq 2$. Let
$\overline{\mathbb{C}(z)}$ be the algebraic closure of
$\mathbb{C}(z)$. For each $i=1,...,m,$ select a function $w_{i}\in
\overline{\mathbb{C}(z)}$ with

\

\

$w_{i}^{N}=\prod_{j=1}^{s}(z-z_{j})^{\widetilde{r}_{ij}}$ for $i=
1,\cdots, m$\\

in $\mathbb{C}(z)[w_{1},...,w_{m}]$. Note that $z_j\in
\mathbb{C}$.

\

\

Here $\widetilde{r}_{ij}$ is the lift of $r_{ij}$ to $\mathbb{Z}
\cap [0,N)$. We impose the condition that the sum of the columns
of $A$ are zero. This implies that the cover is \emph{not}
ramified over infinity. The matrix $A$ will be called the matrix
of the covering. Note that our notations here are mostly that of
[W]. Also we consider the row and column spans of $A$ as modules
over the ring $\mathbb{Z}/N\mathbb{Z}$ and so all of the
operations with rows and columns will be carried out in the ring
$\mathbb{Z}/N\mathbb{Z}$, i.e., it will be considered modulo $N$.
The abelian Galois group $G$ of the covering is isomorphic to the
column span of the matrix $A$ and hence can be considered as a
subgroup of $(\mathbb{Z}/N\mathbb{Z})^{m}$ (denoted also by
$\mathbb{Z}_{N}^{m}$). See \cite{12}, 2.2. Sometimes for
preventing confusions we use the symbol $[]_{N}$ to show in which
ring we are working and for example write $[r_{ij}]_{N}$ instead
of $r_{ij}$. The abelian cover is the smooth projective algebraic
curve with function field $\mathbb{C}(z)[w_{1},...,w_{m}]$. It can
easily be seen that any Galois cover of $\mathbb{P}^{1}$ with
abelian Galois (or deck) group is obtained in this way from a
certain matrix $A$. The local monodromy about the branch point
$z_{j}$ is given by the column vector $(r_{1j},....,r_{mj})^{t}$
and thus the order of ramification over $z_{j}$ is $
\frac{N}{gcd(N,\widetilde{r}_{1j},..,\widetilde{r}_{mj})}$. Using
this and the Riemann-Hurwitz formula, the genus $g$ of the cover
is then given by:

\

\

$g= 1+ d(\frac{s-2}{2}- \frac{1}{2N}\sum_{j=1}^{s}
gcd(N,\widetilde{r}_{1j},...,\widetilde{r}_{mj}))$

\

\

Where $d$ is the degree of the covering. Note that the degree $d$
of the covering (or equivalently the order of the Galois group)
can be realized as the size of the row span (equivalently column
span) of the matrix $A$. See \cite{12}.

\

\

\textbf{Remark 2.1.} Note that a $G$-Galois cover $f: C\rightarrow
\mathbb{P}^{1}$ branched above the points $S=\{z_{1},...,z_{s}\}$
corresponds to a surjection $\phi : \pi_{1}(\mathbb{P}^{1}
\setminus S)\twoheadrightarrow G $ (See \cite{96}, Theorem 5.14).
The fundamental group $\pi_{1}(\mathbb{P}^{1} \setminus S)$ is
generated by the loops $\gamma_{j}$ around $z_{j}$ with only the
condition that $\gamma_{1}...\gamma_{s}=1$ and the local monodromy
around $z_{j}$ is given by $\phi(\gamma_{j})$. Our assumption in
the above that the columns sum is zero implies simply that $\infty
\notin S$.

\

\

Next we turn to define families of curves which are abelian
coverings of the projective line.

Let $T \subset (\mathbb{A}^{1})^{s}$ be the complement of the big
diagonals, i.e., $T=\mathcal{P}_{s}= \{(z_{1},....,z_{s})\in
(\mathbb{A}^{1})^{s}\mid z_{i}\neq z_{j} \forall i\neq j \}$. Over
this open affine set, we define a family of abelian covers of
$\mathbb{P}^{1}$ to have the equation:

\

\

$w_{i}^{N}=\prod_{j=1}^{s}(z-z_{j})^{\widetilde{r}_{ij}}$ for $i=
1,\cdots, m$.

\

Where the tuple $(z_{1},...,z_{s})\in T$ and $\widetilde{r}_{ij}$
is the lift of $r_{ij}$ to $\mathbb{Z}\cap [0,N)$ as before. In
this way each $w_{i}$ defines a cyclic cover of $\mathbb{P}^{1}$.
Note that we have presented the abelian cover as a fiber product,
over $\mathbb{P}^{1}$, of cyclic covers given by the above
formulas. Varying the branch points we get a family $f:C\to T$ of
smooth projective curves over $T$ whose fibers $C_t$ are abelian
covers of $\mathbb{P}^{1}$ introduced above.

\

If $f:C\rightarrow T$ is a family of abelian covers constructed as
above, we write $J\rightarrow T$ for the relative Jacobian of $C$
over $T$. This family gives a natural map $j:T\rightarrow A_{g}$.
Let $Z= Z(N,s,A)$ be the closure $\overline{j(T)}$ in $A_{g}$.
Such a family therefore gives rise to a closed subvariety $Z=
Z(N,s,A)$ in the moduli space $A_{g}$ and we have $\dim Z=s-3$.
This holds because any three points on $\mathbb{P}^{1}$ can be
moved to the points $0,1, \infty$. We call the subvariety $Z$ the
\emph{moduli variety} associated to the family $f:C\rightarrow T$.

\

\

\textbf{Remark 2.2.} Consider two families of abelian covers with
matrices $A$ and $A^{\prime}$ over the same
$\mathbb{Z}/N\mathbb{Z}$. If $A$ and $A^{\prime}$ have equal row
spans then the two families are isomorphic. For more details see
\cite{12}.

\

\subsection{The local system associated to an abelian cover}

\

\

In this section we are going to describe an alternative
construction of abelian coverings using line bundles and local
systems. This construction resembles that of \cite{07} in the case
of cyclic coverings of algebraic varieties. Let $G$ be a finite
abelian group. We denote by $\mu_{G}$ the group of the characters
of $G$, i.e., $\mu_{G}=Hom(G,\mathbb{C}^{*})$. Consider a Galois
covering $\pi: X\rightarrow \mathbb{P}^{1}$ with Galois group $G$.
The group $G$ acts on the sheaves $\pi_{*}(\mathcal{O})$ and
$\pi_{*}(\mathbb{C})$ via its characters and we get corresponding
eigenspace decompositions $\pi_{*}(\mathcal{O})=\oplus_{\chi \in
    \mu_{G}} \pi_{*}(\mathcal{O})_{\chi}$ and
$\pi_{*}(\mathbb{C})=\oplus_{\chi \in \mu_{G}}
\pi_{*}(\mathbb{C})_{\chi}$. Let
$L^{-1}_{\chi}=\pi_{*}(\mathcal{O}_{X})_{\chi}$ and
$\mathbb{C}_{\chi}= \pi_{*}(\mathbb{C})_{\chi}$ denote the
eigensheaves corresponding to the character $\chi$. $L_{\chi}$ is
a line bundle and outside of the branch locus of $\pi$,
$\mathbb{C}_{\chi}$ is a local system of rank 1. Following
\cite{91}, we call the bundles $L_{\chi}$ and $z_{j}, j=1,..,s$,
considered as divisors in $\mathbb{P}^{1}$ the \emph{building
data} of the cover. The reason is that these data determine the
cover completely.

\

\

\textbf{Remark 2.1.1.} Note that if $G$ is a finite abelian group,
then $\mu_G= Hom(G,\mathbb{C}^{*})$ is isomorphic to $G$. To see
this, first assume that $G=\mathbb{Z}/N$ is a cyclic group. We fix
an isomorphism between $\mathbb{Z}/N$ and the group of $N$-th
roots of unity in $\mathbb{C}^{*}$ via $1\mapsto exp(2\pi i/N)$.
Now the group $\mu_G$ is isomorphic to this latter group via $\chi
\mapsto \chi(1)$. In the general case, we can extend this to an
isomorphism $\varphi_G: G \cong \mu_G$ because $G$ is a product of finite cyclic groups.\\

For our applications, with notations as in the previous pages, we
fix an isomorphism of $G$ with a product of $\mathbb{Z}/n$'s and
an embedding of $G$ into $(\mathbb{Z}/N)^{m}$.\\

Let $l_{j}$ be the $j$-th column of the matrix $A$. As mentioned
earlier, the group $G$ can be realized as the column span of the
matrix $A$. Therefore we may assume that $l_{j}\in G$. For a
character $\chi$, $\chi(l_{j})\in \mathbb{C}^{*}$ and since $G$ is
finite $\chi(l_{j})$ will be a root of unity. Let
$\chi(l_{j})=e^{\frac{2\alpha_{j}\pi i}{N}}$, where $\alpha_{j}$
is the unique integer in $[0,N)$ with this property. Equivalently,
the $\alpha_{j}$ can be obtained in the following way: let $a\in
G\subseteq (\mathbb{Z}/N)^{m}$ be the element corresponding to
$\chi$ under the above isomorphism. We regard $a$ as an $1\times
m$ matrix. Then the matrix product $a.A$ is meaningful and
$a.A=(\alpha_{1},...,\alpha_{s})$. Here all of the operations are
carried out in $\mathbb{Z}/N$ but the $\alpha_{j}$ are regarded as
integers in $[0,N)$. Using the above facts, we occasionally
consider a character of $G$ as an element
of this group without referring to isomorphism $\varphi_{G}$.\\

Let us denote by $\omega_X$ the canonical sheaf of $X$. Similar to
the case of $\pi_{*}(\mathcal{O}_X)$, the sheaf
$\pi_{*}(\omega_X)_{\chi}$ decomposes according to the action of
$G$. For the line bundles $L_{\chi}$ corresponding to the
character $\chi$ associated to the element $a\in G$ and
$\pi_{*}(\omega_X)_{\chi}$ we have:

\

\

\textbf{Lemma 2.1.2.} With notations as above $L_{\chi}=
\mathcal{O}_{\mathbb{P}^{1}}(\sum_{1}^{s}
<\frac{\alpha_{j}}{N}>)$.\\

where $<x>$ denotes the fractional part of the real number $x$.

\

and \\

$\pi_{*}(\omega_X)_{\chi}= \omega_{\mathbb{P}^{1}} \otimes
L_{\chi^{-1}}=\mathcal{O}_{\mathbb{P}^{1}}(-2+\sum_{1}^{s}<-\frac{\alpha_{j}}{N}>)$.

\

\textbf{Proof.} Note that since the sum of the columns of the
matrix $A$ is zero, the above sum is an integer. One can easily
see that each section of the line bundle
$\mathcal{O}_{\mathbb{P}^{1}}(\sum_{1}^{s}
<\frac{\alpha_{i}}{N}>)$ is a function on which the Galois group
acts as $\chi$ and conversely any such section must be a function
of the above form. The rest of the lemma is \cite{98}, Proposition
1.2.  \hspace{1.5 cm} $\Box$

\

\textbf{Remark 2.1.3.} Note that in the case of a cyclic cover of
$\mathbb{P}^{1}$ the bundles $L_{\chi}$ and local systems
$\mathbb{C}_{\chi}$ coincide with the bundles $\mathcal{L}^{(j)}$
and local systems $\mathbb{L}_{j}$ in \cite{09} as one expects
naturally.

\

\

\subsection{ Two Shimura subvarieties containing $Z$}

\

\

In this subsection we describe two naturally constructed Shimura
subvarieties of $A_{g}$ associated to a family $f:C\rightarrow T$
of abelian covers that contain the moduli variety $Z$. The
construction that we describe below and Lemma 2.2.1 are based on
the works in \cite{11} and \cite{10} which we just reformulate for
abelian covers. For the sake of completeness, we give the proof of
Lemma 2.2.1 in detail.

\

Let us first sketch the construction of $A_g$, the moduli space of
principally polarized abelian varieties of dimension $g$, as a
Shimura variety. Let $V_{\mathbb{Z}}:=\mathbb{Z}^{2g}\subset
V_{\mathbb{Q}}:=\mathbb{Q}^{2g}$ and let $\psi:
V_{\mathbb{Z}}\times V_{\mathbb{Z}}\to \mathbb{Z}$ be the standard
symplectic form. Let $L=Gsp(V_{\mathbb{Z}},\psi)$ be the group of
symplectic similitudes and
$\mathbb{S}:=Res_{\mathbb{C}/\mathbb{R}}\mathbb{G}_m$ be the
Deligne torus. Consider the space $\mathcal{H}_g$ of homomorphisms
$h:\mathbb{S}\to L_{\mathbb{R}}$ that define a Hodge structure of
type $(1,0)+(0,1)$ on $V_{\mathbb{Z}}$ with $\pm (2\pi i)\psi$ as
a polarization. The pair $(L_{\mathbb{Q}},\mathcal{H}_g)$ is a
Shimura datum and $A_g$ can be described as the Shimura variety
associated to this Shimura datum as follows. Let $K_n:=\{g\in
G(\widehat{\mathbb{Z}})|g\equiv 1$ $(mod$ $n)\}$ with $n\geq 3$.
Then $A_{g,n}(\mathbb{C})=L(\mathbb{Q})\setminus
\mathcal{H}_g\times L(\mathbb{A}_f)/K_n$. The natural number $n$
is called the level structure.
Since it does not play any role in the sequel we drop it from the notation.\\

As $A_g$ has the structure of a Shimura variety, one can talk
about its Shimura (or special) subvarieties. Consider an algebraic
subgroup $N\subset L_{\mathbb{Q}}$ for which the set \\

$Y_N=\{h\in \mathcal{H}_g|$ $h$ factors through $N_{\mathbb{R}}\}$\\

is non-empty. If $Y^{+}$ is a connected component of $Y_N$ and
$\gamma K_n\in L(\mathbb{A}_f)/K_n$, the image of $Y^{+}\times
\{\gamma K_n\}$ in $A_g$ is an algebraic subvariety. We define a
\emph{Shimura subvariety} as an algebraic subvariety $S$ of $A_g$
which arises in this way, i.e., there exists a connected component
$Y^+\subset Y_N$ and an element $\gamma K_n\in
L(\mathbb{A}_f)/K_n$
such that $S$ is the image of $Y^{+}\times \{\gamma K_n\}$ in $A_g$.\\

Now consider a family $C\to T$ of abelian covers of
$\mathbb{P}^{1}$ with abelian Galois group $G$ as in the beginning
of section $2$. The Jacobians in the family $J\rightarrow T$ admit
naturally an action of the group ring $\mathbb{Z}[G]$. This action
defines a Shimura subvariety of PEL type $S(G)$ in $A_{g}$ that
contains $Z$. More precisely, fix a base point $t \in T$ and let
$(J_t,\lambda)$ be the corresponding Jacobian with $\lambda$ as
its principal polarization. Let $(V_{\mathbb{Z}},\psi)$ be as
above. We fix a symplectic similitude $s :H^{1}(J_t,
\mathbb{Z})\to V_{\mathbb{Z}}$. Let $F=\mathbb{Q}[G]$. Via $s$,
the Hodge structure on $H^1(J_t,\mathbb{Q})$ corresponds to a
point $y\in \mathcal{H}_g$ and one obtains the structure of an
$F$-module on $V_{\mathbb{Q}}$. $F$ is isomorphic to a product of
cyclotomic fields and is equipped with a natural involution $*$
which is complex conjugation on each factor. The polarization
$\psi$ on $V_{\mathbb{Q}}$ satisfies

\

$\psi(bu,v)=\psi(u,b^{*}v)$ for all $b\in F$ and $u,v \in V$.

\

Define the subgroup $N$ as in \cite{10}:

\

$N= Gsp(V_{\mathbb{Q}}, \psi)\cap GL_{F}(V_{\mathbb{Q}})$

\

If $h_{0}: \mathbb{S}\rightarrow Gsp_{2g, \mathbb{R}}$ is the
Hodge structure on $V=H^{1}(J_{t},\mathbb{Z})$ corresponding to
the point $y \in \mathcal{H}_{g}$, then by the above $F$-action
this homomorphism factors through the subvariety $N_{\mathbb{R}}$.
Define the set $Y_N$ as above. The point $y$ lies in $Y_{N}$ and
there is a connected component $Y^{+}\subseteq Y_{N}$ which
contains $y$. We define $S(G)$ to be the image of
$Y^{+}$ under the map\\
$\mathcal{H}_g\to L(\mathbb{Z})\setminus \mathcal{H}_g\cong
L(\mathbb{Q})\setminus \mathcal{H}_g\times
L(\mathbb{A}_f)/L(\widehat{\mathbb{Z}})\cong A_g(\mathbb{C})$

\

\

Since $Z \subseteq S(G)$, we have that $s-3 \leq \dim S(G)$.
Therefore if $\dim S(G)=s-3$, it follows that $Z=S(G)$ and hence
$Z$ will be a Shimura subvariety of $A_{g}$. We can find the
dimension of the variety $S(G)$ by finding the tangent space to it
at an arbitrary point. To do this we have to consider the
eigenspaces of the action of the group $G$ on cohomology. The
group $G$ acts on the cohomology $H^{1}(C_{t},\mathbb{C})$ by its
characters. There is also a natural action on
$H^{1,0}=H^{0}(C_{t}, \omega_{C_{t}})$. Let $H^{0}(C_{t},
\omega_{C_{t}})_{(n)}$ be the subspace on which the group $G$ acts
by character $n\in \mu_{G}$. By Remark 2.1.1, there is a
corresponding element in $G$ which we again denote by $n$. Put
$d_{n}= \dim_{\mathbb{C}}H^{0}(C_{t}, \omega_{C_{t}})_{(n)}$. We
have:

\

\

\textbf{Lemma 2.2.1} $\dim S(G)= \sum_{2n\neq 0} d_{n}d_{-n}+
\frac{1}{2}\sum_{2n=0}d_{n}(d_{n}+1)$.

\

Note that $2.0=0$ in $G$ and $d_{0}=0$, so in fact the second sum
in the right hand side of the above equality is always meaningful
and if $|G|$ is an odd number it will be zero.

\

\textbf{Proof.} We calculate $\dim T_{y}(Y_{N})$ at the point $y
\in \mathcal{H}_{g}$. The dimension of the tangent space of $S(G)$
at the point $y$ will be equal to this number. To compute $\dim
T_{y}(\mathcal{H}_{g})$, we note that:

\

$T_{y}(\mathcal{H}_{g})= Hom^{sym}(H^{1,0},
V_{\mathbb{C}}/H^{1,0}):=$

\

$\{\beta : H^{1,0}\rightarrow V_{\mathbb{C}}/H^{1,0} \mid
\overline{\phi}(v,\beta(v^{\prime}))=\overline{\phi}(v^{\prime},\beta(v))
\forall v,v^{\prime} \in H^{1,0}\}$

\

i.e., the elements of $T_{y}(\mathcal{H}_{g})$ that are their own
dual via the isomorphisms induced by $\overline{\phi}$.  For a
more detailed discussion see \cite{10}. The subspace $T_{y}(Y_{N})
\subset T_{y}(\mathcal{H}_{g})$ consists of $\beta \in
Hom^{sym}(H^{1,0}, V_{\mathbb{C}}/H^{1,0})$ (symmetric with
respect to $\overline{\phi}$) that respect the $F$-action on $V$,
that is, are $F_{\mathbb{C}}$-linear. Any such $\beta$ can be
written as the sum $\sum \beta_{n}$, where
$\beta_{n}:H^{1,0}_{\mathbb{C},n} \rightarrow H^{0,1}_{\mathbb{C},
    n}$ is the induced action on the eigenspaces. These $\beta_{n}$
should satisfy the relation

\

$\overline{\phi}_{n}(v, \beta_{-n}(v^{\prime}))=
\overline{\phi}_{-n}(v^{\prime}, \beta_{n}(v)) $.

\

Note that the map $\overline{\phi}_{n}$ induced by the
polarization $\phi$ gives a duality between
$H^{1,0}_{\mathbb{C},n}$ and $H^{0,1}_{\mathbb{C}, (-n)}$. So we
have a duality between $\beta_{n}$ and $\beta_{-n}$ if $n\neq -n$
in $G$. If $n= -n$ in $G$, i.e., if $2n=0$ in $G$ this gives a
self duality for $\beta_{k}$. Therefore $\dim T_{y}(Y_{N})$ is
equal to:

\

$ \sum_{2n\neq 0} d_{n}d_{-n}+
\frac{1}{2}\sum_{2n=0}d_{n}(d_{n}+1)$  \hspace{1.5cm}$\Box$

\

\

\textbf{Construction 2.2.2.} The construction of the second
Shimura subvariety that contains $Z$  is in fact Mumford's
construction of "variety of Hodge type ". Namely, let $M$ be the
generic Mumford-Tate group of the family $f:C\rightarrow T$. For
the definition and construction of the generic Mumford-Tate group
we refer to \cite{15} or \cite{09}. Note that $M$ is a reductive
$\mathbb{Q}$-algebraic group . Let $S_{f}$ be the natural Shimura
variety associated to $M$. Note that in general $M\subseteq N$,
where $N$ is the group defined before Lemma 2.2.1, and we have
$S_{f}\subseteq S(G)$. The Shimura subvariety $S_{f}$ is in fact
the \emph{smallest} Shimura subvariety that contains $Z$ and its
dimension depends on the real adjoint group $M^{ad}_{\mathbb{R}}$.
Namely, if $M^{ad}_{\mathbb{R}}= Q_{1}\times...\times Q_{r}$ is
the decomposition of $M^{ad}_{\mathbb{R}}$ to $\mathbb{R}$-simple
groups then $\dim S_{f}= \sum \delta(Q_{i})$. If
$Q_{i}(\mathbb{R})$ is not compact then $\delta(Q_{i})$ is the
dimension of the corresponding symmetric space associated to the
real group $Q_{i}$ which can be read from Table V of \cite{78}. If
$Q_{i}(\mathbb{R})$ is compact, i.e., if $Q_{i}$ is anisotropic we
set $\delta(Q_{i})=0$. We remark that for $Q=PSU(p,q)$,
$\delta(Q)=pq$ and for $Q=Psp_{2p}$, $\delta(Q)=\frac
{p(p+1)}{2}$. Note also that $Z$ is a Shimura subvariety if and
only if $\sum\delta(Q_{i})=s-3$, i.e., if and only if $\dim Z=\dim
S_{f}=s-3$.

\

\

\textbf{Computation of $d_{n}$.} We have seen that the dimension
of the Shimura variety $S(G)$ can be expressed in terms of the
dimension of the eigenspaces of Galois action on the cohomology of
the fibers. We will now try to compute these dimensions. Let $\chi
\in\mu_{G}$ be a character of $G$ and $n=(a_{1},..,a_{m}) \in G$
be the element corresponding to $\chi$ under the isomorphism
$\varphi_{G}$ introduced in Remark 2.1.1. Let
$\alpha_{1}$,...,$\alpha_{s}$ be as in Lemma 2.1.2. We have the
following result:

\

\

\textbf{Proposition 2.2.3.} For $C$ an abelian cover of
$\mathbb{P}^{1}$ we have:

\

\

$d_{n}=h^{1,0}_{\chi}(C)=-1+\sum_{1}^{s}<- \frac{\alpha_{j}}{N}>$

\

\textbf{Proof.} By Lemma $2.1.2$ we have that:

\

$\pi_{*}(\omega_C)_{\chi}= \omega_{\mathbb{P}^{1}} \otimes
L_{\chi^{-1}}=\mathcal{O}_{\mathbb{P}^{1}}(-2+\sum_{1}^{s}<-
\frac{\alpha_{j}}{N}>)$

\

It follows that:

\

$h^{1,0}_{\chi}(C)= h^{0}(\pi_{*}(\omega_C)_{\chi})=
-1+\sum_{1}^{s}<- \frac{\alpha_{j}}{N}>$. \hspace{1.5cm}$\Box$

\

\

Note that it naturally follows that

\

$d_{-n}=h^{1,0}_{\chi^{-1}}(C)=h^{0,1}_{\chi}(C)=-1+\sum_{1}^{s}<\frac{\alpha_{j}}{N}>$.

\

\

\textbf{Remark 2.2.4.} There are other methods to compute the
dimension of the eigenspaces $d_{n}$. Note that the abelian Galois
group $G$ of the covering is a (possibly proper) subgroup of
$\mathbb{Z}^{m}_{N}$ and therefore we can show an element of $G$
as an $m$-tuple $n=(n_{1},...,n_{m})$. The space $H^{1,0}_{n}$ of
differential forms on which $G$ acts via character $n$ is
generated over $\mathbb{C}(z)$ by the differential form
$w_n=\prod(z-z_{j})^{-t_{j}(-n)}dz$, where
$t_{j}(n)=<\frac{\sum_{i=1}^{m} n_{i}\widetilde{r}_{ij}}{N}> $ and
$<.>$ denotes the fractional part of a real number. It is then
straightforward to check that a meromorphic form
$p(z)\prod(z-z_{j})^{-t_{j}(-n)}dz$ is holomorphic if and only if
$p(z)$ is a polynomial of degree at most $t(n)= \sum t_{j}(-n)-1$,
see \cite{12}, Lemma 2.6. So the dimension of $H^{1,0}_{n}$ is
equal to $t(n)$. Alternatively, one can use the Chevalley-Weil
formula to compute the dimension of the eigenspaces. See
\cite{34}.

\

\

\section{Examples of Shimura varieties arising from abelian covers}

In \cite{11}, Moonen completed the list of Shimura subvarieties
generated by families of cyclic covers of $\mathbb{P}^{1}$ and
proved that in the locus of cyclic covers of $\mathbb{P}^{1}$
there are no more Shimura varieties. The fiber genus of these
families is bounded by $8$ confirming the bound given by the
corrected version of Coleman conjecture. See page 2. In \cite{10},
Oort and Moonen give a table of $7$ examples of abelian non-cyclic
Galois covers of $\mathbb{P}^{1}$ that generate Shimura
subvarieties in $A_{g}$. All of these examples satisfy the
equality $\dim S(G)=s-3$. Their method to obtain these examples is
based on analyzing the decomposition of Jacobians up to the
isogeny under the action of group  ring $\mathbb{Q}[G]$. Our
argument here is based on checking the equality of Lemma 2.2.1.
Checking whether this equality holds is something that can be
checked on a computer and by using a computer program we have
checked the examples which satisfy this equality. However, our
computer search did not provide a further example satisfying $\dim
S(G)=s-3$. We are able however to prove that for $s=4$ the table
contains all examples with $\dim S(G)=s-3=1$ (see Theorem 3.1
below). We moreover study the families that \emph{do not} appear
in the our table, i.e., those that do not satisfy the equality
$\dim S(G)=s-3$. In this case, $Z\neq S(G)$ but it does not imply
that $Z$ is not a Shimura subvariety: it could still be a smaller
Shimura subvariety (inside $S(G)$) or in other words, there might
be Hodge classes that are \emph{not} given by the action of
$\mathbb{Z}[G]$. We are able to show that some large classes of
families, including all families with $s=4$, do not give rise to
Shimura subvarieties in $A_{g}$ provided that the family is
irreducible, see Theorem 6.1.1.

\

\

\begin{table}

    \

    \

    \caption{Monodromy data of families of abelian coverings that
        generate Shimura subvarieties. }

    \

    \

    \begin{tabular}{|c|c|c|c|}

        \hline
        genus & Galois group& N & monodromy data\\

        \hline
        1& $\mathbb{Z}/2\mathbb{Z} \times \mathbb{Z}/2\mathbb{Z} $ & 2& $\{(1,0)(1,0)(0,1)(0,1)\}$\\
        \hline
        2 & $\mathbb{Z}/2\mathbb{Z} \times \mathbb{Z}/2\mathbb{Z} $ & 2& $\{(1,0)(1,0)(1,0)(1,1)(0,1)\}$\\
        \hline
        3 &  $\mathbb{Z}/2\mathbb{Z} \times \mathbb{Z}/4\mathbb{Z} $& 4& \{(2,0)(2,1)(0,1)(0,2)\}\\

        \hline

        3 &  $\mathbb{Z}/2\mathbb{Z} \times \mathbb{Z}/4\mathbb{Z} $& 4& \{(2,0)(2,2)(0,1)(0,1)\}\\

        \hline
        3 &  $\mathbb{Z}/2\mathbb{Z} \times \mathbb{Z}/2\mathbb{Z} $ & 2& \{(1,0)(1,0)(1,1)(1,1)(0,1)(0,1)\}\\

        \hline
        4 & $\mathbb{Z}/2\mathbb{Z} \times \mathbb{Z}/6\mathbb{Z} $ & 6& \{(3,0)(3,1)(0,2)(0,3)\} \\

        \hline
        4 & $\mathbb{Z}/3\mathbb{Z} \times \mathbb{Z}/3\mathbb{Z} $  & 3& \{(1,0)(1,0)(1,2)(0,1)\}\\
        \hline

    \end{tabular}
\end{table}

\

\

\textbf{Theorem 3.1.} The families in Table $1$ give rise to
Shimura subvarieties in $A_{g}$.  Moreover for $s=4$ this table
contains all examples of families of abelian non-cyclic Galois
covers of $\mathbb{P}^{1}$ for which $\dim S(G)=1 (=s-3)$ and
contains all such examples for $5\leq s \leq 7$, $2\leq m \leq 5$
and $N\leq 20$.

\

\

\textbf{Proof.} One can compute the dimensions $d_{n}$ of
eigenspaces with the aid of the formula in Proposition 2.2.3. It
is straightforward to check that in all of these cases $\dim
S(G)=\dim Z=s-3$ and therefore $Z=S(G)$ is a Shimura subvariety of
$A_{g}$. By using the computer program mentioned on the last page,
one can check that for the triples $(m,s,N)$ satisfying the bounds
of the theorem the only possibilities are the ones listed above.
For $s=4$, note that if $\dim S(G)=1$ then by using Lemma 2.2.1,
there is a unique $a\in G$ such that $d_{a}=d_{-a}=1$ and for all
other $n\in G$, $d_{n}d_{-n}=0$. We may therefore assume that the
first row of the matrix $A$ satisfies this equality. By results of
\cite{11} or \cite{09} we know that there are only finitely many
of these with $N\leq 12$ and by the aforementioned computer
program we may check that the above examples are the only ones
which satisfy $\dim S(G)=1$ and $N\leq 12$. This means that Table
1 contains all examples with $\dim S(G)=1$. By using the same
computer program we see that this table contains all examples with
$\dim S(G)=s-3$ and with conditions as above. \hspace{1.5cm}$\Box$

\

\

\section{The Dwork-Ogus obstruction}

As we remarked earlier, we are going to exclude further examples
of Shimura subvarieties arising from families of abelian covers.
To do this we will need an obstruction introduced by Dwork and
Ogus in \cite{86}. In this section we follow \cite{11}, \S 5 which
involves the same ideas and techniques as in \cite{89}, \S 5 using
the Dwork-Ogus obstruction. Although we will encounter cases that
we can not use this obstruction, it remains a crucial tool in
proving that a certain variety is not a Shimura subvariety. The
construction of the obstruction is as follows:

\

\

Let the $f:C\rightarrow T$ be a family of smooth projective curves
with an irreducible base scheme $T$. We denote the sheaf of
relative differentials with $\omega_{C/T}$ and the Hodge bundle
$\mathbb{E}=\mathbb{E}(C/T)=f_{*}\omega_{C/T}$. Consider the
Kodaira-Spencer map $\kappa: Sym^{2}(\mathbb{E})\rightarrow
\Omega^{1}_{T}$ (usually the dual of this map is defined to be the
Kodaira-spencer map, but as we will mainly work with this map,
rather than the original Kodaira-Spencer, we name it as such). The
multiplication map $mult:Sym^{2}(\mathbb{E})\rightarrow
f_{*}(\omega^{\otimes2}_{C/T})$ induces the following sheaves:

 \

 \

 $\mathcal{K}=Ker(mult)=Ker(Sym^{2}(\mathbb{E})\rightarrow f_{*}(\omega^{\otimes2}_{C/T}))  $

 \

 $\mathcal{L}= Coker(mult^{\vee})=Coker((f_{*}\omega^{\otimes2}_{C/T})^{\vee}\rightarrow
 Sym^{2}(\mathbb{E})^{\vee})$,

 \

 If the fibers are not hyperelliptic, by a famous result of Max
 Noether (see \cite{80}), $mult$ is surjective and $\mathcal{K}$ is dual to $\mathcal{L}$.

 \

 Now, if $C$ is a smooth projective curve over a field $k$ of
 characteristic $p>0$ with an ordinary Jacobian (we call such a curve an\emph{ ordinary curve})
 and a principal polarization $\lambda$, Serre-Tate theory (see \cite{86})
 guaranties that there exists a canonical lifting $J^{can}$ of $J$ to the Witt
 ring $W(k)$. The question of whether the canonical lifting of a
 Jacobian $J$ is again a Jacobian has been of main interest
 and Dwork and Ogus have shown in \cite{86} that even over the Witt ring
 of length $2$ this is a very restrictive condition and in
 general is not true. Their method consists of constructing an
 obstruction $\beta$ such that $\beta=0$ if and only if the
 canonical lifting $J^{can}$ is a Jacobian. They then show that
 this obstruction is generically non-zero. We recall the
 construction of $\beta$ in short. The curve $C$ is called
 pre-$W_{2}-canonical$ if the canonical lifting
 $(J^{can},\lambda^{can})$ over $W_{2}(k)$ is isomorphic to the
 Jacobian of a smooth projective curve $Y$ as a principally
 polarized abelian variety. According to Dwork-Ogus theory, the
 obstruction $\beta_{C}$ to the existence of such $Y$ is the restriction of an element
 $\beta_{C}\in Sym^{2}(\mathbb{E})^{\vee}$ to the kernel
 $ker(mult)$ of the multiplication map. This obstruction can be
 generalized to an obstruction for families $f:C\rightarrow T $ of
 ordinary curves to give an obstruction $\widetilde{\beta}_{C/T}$
 which is a global section of $F^{*}_{T}\mathcal{L}(C/T)$ where
 $F_{T}:T\rightarrow T$ denotes the absolute Frobenius map and the value of
 $\widetilde{\beta}_{C/T}$ at $t\in T$ is equal to $F^{*}_{k}(\beta_{C_{t}/k})$. Note
 that since the family is assumed to have ordinary fibers the
 inverse Cartier operator $\gamma :F^{*}_{T}\mathbb{E}\rightarrow
 \mathbb{E}$ is an $\mathcal{O}_{T}$-linear isomorphism. The matrix of the inverse transpose of $\gamma$, which
 is the Frobenius action on $R^{1}f_{*}\mathcal{O}_{C}$, is called the \emph{Hasse-Witt matrix} of the family. By a result
 of Katz in \cite{71}, the pull-back $F^{*}_{T}\mathcal{L}(C/T)$ comes
 equipped with a natural flat connection:

 \

 $\nabla: F^{*}_{T}\mathcal{L}\rightarrow F^{*}_{T}\mathcal{L}\otimes
 \Omega^{1}_{T/k}$.

 \

 For the Dwork-Ogus obstruction
 $\widetilde{\beta}_{C/T} $, it holds that $-\nabla
 \widetilde{\beta}_{C/T}: F^{*}_{T}\mathcal{K} \rightarrow
 \Omega^{1}_{T/k}$ is equal to the composition

 \

$F^{*}_{T}\mathcal{K}\hookrightarrow
F^{*}_{T}Sym^{2}(\mathbb{E})\xrightarrow{Sym^{2}(\gamma)}
Sym^{2}(\mathbb{E})\xrightarrow{\kappa} \Omega^{1}_{T/k}$ (*)

 \

  \

 The map $\kappa: Sym^{2}(\mathbb{E})\rightarrow\Omega^{1}_{T/k}$
 is the Kodaira-Spencer map associated to the family
 $f:C\rightarrow T$.

 \

 The above description of $\nabla \widetilde{\beta}_{C/T}$ gives us
 something computable which we will use later to show that the
 obstruction is not zero for our families.

 \

 The curves and their families which we introduced in section 1
 were only defined over $\mathbb{C}$. In order to be able to use
 the above characteristic $p$ tools, we need to set the scene in
 such a way that this reduction mod $p$ makes sense.\\

 Let $f:C\rightarrow T$ be a family of abelian covers as in section
 1. Let $R=\mathbb{Z}[1/N,u]/\Phi_N$, where $\Phi_N$ is the $N$th
 cyclotomic polynomial. Note that $R$ can be embedded into
 $\mathbb{C}$ by sending the image of $u$ to $\exp (2\pi i/N)$. We
 consider $T\subset (\mathbb{A}^1_R)^s$ as the complement of the
 big diagonals, i.e., as the $R$-scheme of ordered $s$-tuples of
 distinct points in $\mathbb{A}^1_R$. For a prime number $p$, we
 denote by $\wp$ a prime of $R$ lying above $p$. One can choose a
 prime number $p\equiv 1$ (mod $N$) and an open subset $U$ of
 $T\otimes \mathbb{F}_{p}\cong T\otimes_R R/\wp$ such that for all
 $t\in U$, the fibers are ordinary curves in characteristic $p$.
 This is possible for example by \cite{01}, Theorem on page 2. For such
 $p$ and $U$, consider the restricted family $C_{U}\rightarrow U$.
 The abelian group $G$ also acts on the sheaves
 $\mathcal{L}(C_{U}/U)$ and gives the eigensheaf decomposition
 $\mathcal{L}(C_{U}/U)=\oplus_{n\in G}\mathcal{L}_{(n)}$. The same
 is true for $\mathbb{E}_{U}=\mathbb{E}(C_{U}/U)$ and
 $\mathcal{K}_{U}=\mathcal{K}(C_{U}/U)$. This in turn gives us the
 decomposition
 $\widetilde{\beta}_{C_{U}/U}=\sum_{n}\widetilde{\beta}_{n}$. Here
 $\widetilde{\beta}_{n}$ is considered as a section of
 $F^{*}_{U}\mathcal{L}_{n}$.

 \

  \

 The main observation here is that if the family gives rise to a
 Shimura subvariety in $A_{g}$ then the Dwork-Ogus obstruction
 vanishes:

 \

 \

 \textbf{Lemma 4.1.} For prime number $p$ and open subset $U$ as
 above, if the family gives rise to a Shimura subvariety
 $Z\subseteq A_{g}$ then for any $t\in U$ we have that the Jacobian
 $J_{t}$ is pre-$W_{2}$-canonical and in particular
 $\widetilde{\beta}_{C_{U}/U}=0$.

 \

 \

 \textbf{Proof.} This follows from results of \cite{92}. In fact if the
 moduli variety $Z$ is a Shimura subvariety and $t\in T$ is an
 ordinary point (i.e., $J(C_t)$ is an ordinary abelian variety)
 then the canonical lifting $J^{can}_{t}$ of $J_{ t}$ is a
 $W(k)$-valued point of $Z$. This means in particular that it is a
 Jacobian and hence $J_{t}$ is pre-$W_{2}$- canonical. By
 Dwork-Ogus theory this forces $\widetilde{\beta}_{C_{U}/U}$ to be
 zero.\hspace{1.5cm}$\Box$

 \

 Now assuming that the fibers of the family are ordinary over $U$
 and the family gives rise to a Shimura subvariety in $A_{g}$, it
 follows from the Lemma 4.1 that $\widetilde{\beta}_{C_{U}/U}=0$
 and hence $\nabla \widetilde{\beta}_{C_{U}/U}=0$. This shows that
 the composition map (*) should vanish identically.

 \

 From now on we just work with the restricted family $C_{U}/U$
 whose fibers are all ordinary instead of $C/T$ and denote it
 simply as $C/U$. Next we remark that the sequence (*) factors
 through the map

 \

 $Sym^{2}(\mathbb{E})\rightarrow
 Sym^{2}(\mathbb{E})_{(0)}\xrightarrow {mult_{(0)}}
 f_{*}(\omega_{C/U}^{\otimes 2})_{(0)}$

 \

 Where by the index $(0)$ we mean the subspace of invariant
 elements under the action of $G$, i.e., the subspace on which $G$
 acts with the trivial character. Note that
 $f_{*}(\omega_{C/U}^{\otimes 2})_{(0)}$ is a locally free sheaf of
 rank $s-3$. The above factorization follows from the general fact
 that the fiber of $Sym^{2}(\mathbb{E})_{(0)}$ at $t$ can be
 identified with (dual of) the space of $G$-equivariant
 deformations of $C_{t}$, i.e., the deformations for which the
 $G$-action also deforms along. Since there is a $G$-action on our
 whole family, the Kodaira-Spencer map should factor through the
 above map. The last map in the above sequence is just
 multiplication of forms.

\

\

\textbf{Proposition 4.2.} With notations as above, the map

\

\

$F^{*}_{U}\mathcal{K}_{(0)}\hookrightarrow
F^{*}_{U}Sym^{2}(\mathbb{E}_{U})_{(0)}\xrightarrow{Sym^{2}(\gamma)}
Sym^{2}(\mathbb{E}_{U})_{(0)}\xrightarrow{mult_{(0)}}
f_{*}(\omega_{C/U}^{\otimes2})_{(0)}$

\

vanishes identically provided that the family gives rise to a
Shimura subvariety in $A_{g}$.

\

\

\textbf{Proof.} This is basically \cite{11}, Prop 5.8 in the more
general setting of abelian covers. Let us first note that the
induced Kodaira-Spencer map
$\kappa_{(0)}:f_{*}(\omega_{C/T}^{\otimes2})_{(0)}\rightarrow
\Omega^{1}_{T} $ is injective. We remark that this map is in fact
the dual of the usual Kodaira-Spencer map $\rho:
\Theta_{T}\rightarrow H^{1}(\Theta_{C/T})$. Therefore its
injectivity means the surjectivity of the Kodaira-Spencer map,
i.e., the versality (or completeness) of our family. Now if
$D/k[\epsilon]$ is a $G$-equivariant first order deformation of
the fiber $C_{t}$ the versality means that $D/k[\epsilon]$ can be
obtained by pull-back from our family. But this holds because in
this case $D/G$ is isomorphic to $\mathbb{P}^{1}_{k[\epsilon]}$
and so as an abelian cover of $\mathbb{P}^{1}$ it can be obtained
by pull-back from our family. Note that in this proof we use the
versality of this family for the problem of deforming a curve with
$G$-structure (and not the versality for a problem of deforming an
abelian variety). If the fibers are non-hyperelliptic the
vanishing of the above map follows directly from the theory of
Dwork-Ogus, see \cite{86}, together with the injectivity of
$\kappa_{(0)}$ discussed above. In fact, according to \cite{86},
the exact sequence (*) being equal to $-\nabla
\widetilde{\beta}_{C/U}$ vanishes identically (Lemma 4.1 above.
Note that this is where we use the hypothesis that the family
gives rise to a Shimura subvariety). Injectivity of $\kappa_{(0)}$
then gives the vanishing of the claimed map. So we may assume that
the fibers are hyperelliptic curves. Namely, with $\iota \in
Aut(C/U)$ being the hyperelliptic isomorphism we conclude from
results of \cite{80} that although the multiplication map
$Sym^{2}(\mathbb{E}_{U})\rightarrow f_{*}(\omega^{\otimes2})$ is
no longer surjective, the induced map $mult_{\iota}:
Sym^{2}(\mathbb{E}_{U})_{\iota} \rightarrow
f_{*}(\omega^{\otimes2})_{\iota}$ on the the sheaves of invariants
of $\iota$  is again surjective. Since our family is contained in
the hyperelliptic locus, this implies that the map $mult_{(0)}$ is
also surjective and this forces $\widetilde{\beta}_{(0)}$ to be an
$\mathcal{O}_{U}$-linear map. If $\widetilde{\beta}_{(0)}$ is zero
of course $\nabla\widetilde{\beta}_{(0)}$ will be also zero. To
complete the proof one can check that $0$-component analogue of
the exact sequence (*) also holds true for $\nabla
\widetilde{\beta_{(0)}}$ . \hspace{1.5cm}$\Box$

\

\

\

\section{The generalization of a lemma}

\

\

For our classification purposes, we will need the generalization
to the abelian case of a lemma in \cite{11} that concerns only
with cyclic coverings. See \cite{01}, Lemma 5.1.i. This lemma
allows us to compute explicitly the Hasse-Witt matrix of an
abelian covering which considering the above constructions will be
needed to compute the obstruction $\widetilde{\beta}_{C/U}$.
Consider an abelian cover with Galois group $G$ and matrix $A$.
Let $a=(a_{1},...,a_{m}) \in G \subseteq \mathbb{Z}_{N}^{m}$ be an
element in the Galois group of the abelian covering (or the
corresponding character $\chi$, see Remark 2.1.1). Consider
$a.A=(\alpha_{1},...,\alpha_{s})$ as in Lemma 2.1.2. Let
$\widetilde{r}_{ij}$ be the lift of $r_{ij}$ as in section $2$ and
$\widetilde{A}=(\widetilde{r}_{ij})$. Let
$a.\widetilde{A}=(\sum_{1}^{m}a_{i}\widetilde{r}_{i1},...,\sum_{1}^{m}a_{i}\widetilde{r}_{is})=(\widetilde{\alpha}_{1},...,\widetilde{\alpha}_{s})$.

\

\

Next take a prime number $p$ such that $p\equiv 1$ (mod $N$) and
let $q= \frac{p-1}{N}$.

\

\

\textbf{Lemma 5.1.} With notation as above, there exists a basis
$(\xi_{a,i})_{a,i}$ of $H^{1}(Y,\mathcal{O}_{Y})$ with respect to
which the $(h_{\nu \iota})$ entry of the Hasse-Witt matrix of the
abelian covering $Y$ is given by the formula:

\

\

$\sum_{\sum
    l_{i}=\Upsilon}\binom{q.[\alpha_{1}]_{N}}{l_{1}}...\binom{q.[\alpha_{s}]_{N}}{l_{s}}z_{1}^{l_{1}}...z_{s}^{l_{s}}$

\

\

Where $\Upsilon= (d_{n}-\iota)(p-1)+(\nu-\iota)$ and
$\binom{a}{b}=\frac{a!}{b!(a-b)!}$.

\

\

\textbf{Proof.} As in \cite{01} take
$U_{1}=\mathbb{P}^{1}-\{\infty\}$ and $U_{2}=\mathbb{P}^{1}-\{0\}$
and $V_{i}=\pi^{-1}U_{i}$ for $i=1,2$. Let

\

$v_{a}=
w_{1}^{a_{1}}...w_{m}^{a_{m}}(z-z_{1})^{-[\frac{\widetilde{\alpha}_{1}}{N}]}...(z-z_{s})^{-[\frac{\widetilde{\alpha}_{s}}{N}]}$

\

 Then we have:

 \

 \

 $\Gamma(V_{1},\mathcal{O}_{V_{1}})= \bigoplus_{a\in G} k[z]v_{a}$

 \

 $\Gamma(V_{2},\mathcal{O}_{V_{2}})= \bigoplus_{a\in G}
 k[z^{-1}]z^{-|a|-1}v_{a} $

 \

 $\Gamma(V_{1}\cap V_{2})= \bigoplus_{a\in G} k[z,z^{-1}]v_{a}$

 \

 \

 Defining $\xi_{a,i}= z^{-i}v_{a}$ for $i=1,...,|a|$, with
 $|a|=\dim H^{1}(Y,\mathcal{O}_{Y})_{a}$, we see that the
 $\xi_{a,i}$'s form a basis for $H^{1}(Y,\mathcal{O}_{Y})_{a}=
 \frac{\Gamma(V_{1}\cap
     V_{2})_{a}}{\Gamma(V_{1})_{a}+\Gamma(V_{2})_{a}}$ and
 $\{\xi_{a,i}|a\in G\setminus \{0\}, 0<i\leq |a|\}$ is a basis for
 $H^{1}(Y,\mathcal{O}_{Y})$. By the above choice of $p$, the
 Hasse-Witt map $F:H^{1}\rightarrow H^{1}$ induces a corresponding
 map $F_{a}:H^{1}_{a}\rightarrow H^{1}_{a}$. This gives us a
 decomposition of the Hasse-Witt map $F$ corresponding to the
 decomposition of $H^{1}(Y,\mathcal{O}_{Y})$ as a sum of
 eigenspaces described above.

 \

 If $B_{a}$ is the matrix of the Hasse-Witt map
 $F_{a}:H^{1}_{a}\rightarrow H^{1}_{a}$ then the $(i,j)$ entry of
 $B_{a}$ is given by the coefficient of $\xi_{a,j}$ in
 $\xi_{a,i}^{p}$. This follows from the fact that the
 $\xi_{a,i}\otimes 1$ determine a local basis for the bundle
 $F_{T}^{*}(R^{1}f_{*}\mathcal{O}_{C})_{a}$ and the Hasse-Witt
 operator $\gamma : F_{T}^{*}(R^{1}f_{*}\mathcal{O}_{C})_{a}
 \rightarrow (R^{1}f_{*}\mathcal{O}_{C})_{a} $ with respect to
 these bases is given by the p-th power endomorphism of
 $\mathcal{O}_{C}$. This is due to the fact that the $p$-linear
 composite map

 \

 $R^{1}f_{*}\mathcal{O}_{C}\rightarrow
 F_{T}^{*}(R^{1}f_{*}\mathcal{O}_{C}) \rightarrow
 (R^{1}f_{*}\mathcal{O}_{C})$

 \

 is induced by the $p$-th power endomorphism of $\mathcal{O}_{C}$
 (cf. \cite{72}, 2.3.4.1.4). Now one sees that the coefficient of
 $\xi_{a,j}$ in this polynomial is as claimed above:

 \

 \

 $\sum_{\sum
     l_{i}=\Upsilon}\binom{q.[-\alpha_{1}]_{N}}{l_{1}}...\binom{q.[-\alpha_{s}]_{N}}{l_{s}}z_{1}^{l_{1}}...z_{s}^{l_{s}}$. \hspace{1.5cm}$\Box$

 \

 \

 \section{Excluding non-Shimura examples}

 \

 \

 At this point we are going to exclude the families of abelian
 covers of the projective line that do not give rise to Shimura
 subvarieties in $A_{g}$. For some technical reasons working with
 families with $4$ branch points is different from families with
 more branch points and it should be noted that this is in some
 sense the most important case as most of the examples of Shimura
 families that have been found in \cite{11} or \cite{10} are obtained from
 families with $4$ branch points.

 \

 \

 \textbf{Remark 6.1.} From now on we will assume that our families
 are families of \emph{irreducible} abelian covers of
 $\mathbb{P}^{1}$. For cyclic covers this implies that the single
 row of the associated matrix is not annihilated by a non-zero
 element of $\mathbb{Z}/N\mathbb{Z}$. More generally, for abelian
 covers of $\mathbb{P}^{1}$ this implies that the rows of the
 associated matrix are linearly independent over
 $\mathbb{Z}/N\mathbb{Z}$.

 \

 \

 \subsection{The case of four branch points}

  \

 \textbf{Theorem 6.1.1.} Let $Y \rightarrow T$ be a family of
 irreducible abelian covers with $s=4$, i.e., with 4 branch points.
 Then the associated subvariety $Z \subseteq A_{g}$ is a Shimura
 curve if and only if $Z=S(G)$, i.e., if and only if it appears in
 Table 1.

 \

 \

 \textbf{Proof}. Clearly the statement holds if $Z=S(G)$ as $S(G)$
 is a Shimura variety of PEL type. Now assume on the contrary that
 $Z \neq S(G)$ but $Z$ is a Shimura subvariety and we will derive a
 contradiction. By Theorem 3.1, the assumption $Z \neq S(G)$
 implies that $\dim S(G) >1$. Note that in this theorem we have
 classified all cases where $\dim Z=\dim S(G)=1$, i.e., the cases
 for which $s=4$ and $Z=S(G)$. Now suppose that $\dim S(G)>1$. By
 Lemma 2.2.1, if $d_n=0$ or $d_{-n}=0$ then $n\in G$ does not
 contribute to $\dim S(G)$. Also, Lemma 2.2.3 implies that when
 $s=4$, $d_n\leq 2$ and $d_n+d_{-n}\leq 2$. So if $d_n=2$, then
 $d_{-n}=0$ which we excluded. These facts together show that there
 are pairs $a, a^{\prime} \in G $ with $a^{\prime}\neq \pm a$ such
 that $d_{a}=d_{-a}=1$ and $d_{a^{\prime}}=d_{-a^{\prime}}=1$.
 Consequently, the spaces $H^{1,0}_n$ for $n\in \{\pm a, \pm
 a^{\prime}\}$ are $1$-dimensional. For $n\in \{\pm a, \pm
 a^{\prime}\}$, let $w_n$ be the generator of $H^{1,0}_n$ as in
 Remark 2.2.4. For these $n$'s therefore, the Hasse-witt matrix
 $A_{n}$ is a polynomial in $\mathbb{F}_{p}[z_{1},..,z_{4}]$. Note
 that by the discussion just before Lemma 4.1, there exist a
 suitable prime number $p$ and an open subset $U$ of $T\otimes
 \mathbb{F}_{p}$ such that all fibers above $U$ are ordinary.
 Therefore the Hasse-Witt operator is an isomorphism over $U$ and
 so $A_{n}$ is invertible as a section of $\mathcal{O}_{U}$. It
 follows from the description of $w_n$ in Remark 2.2.4 that
 $\omega_{a}.\omega_{-a}=\omega_{a^{\prime}}.\omega_{-a^{\prime}}$
 as a section of the bundle $f_{*}(\omega^{\otimes 2})$. It is a
 non-zero section of the bundle $f_{*}(\omega^{\otimes 2})$ and so
 we must have : \

 \

 $A_{a}.A_{-a}=A_{a^{\prime}}.A_{-a^{\prime}}$

 \

 as polynomials.

 \

 \

 We will show that this identity can not happen with the above
 conditions. The polynomials $A_{n}$ are given by the dual version
 of Lemma $5.1$ and we set $B_{n}=A_{n}\mid_{z_{1}=0}$. It means
 that we have : \

 \

 $B_{n}= \sum_{j_{2}+j_{3}+j_{4}=p-1}
 \binom{q.[-\alpha_{2}]_{N}}{j_{2}}...\binom{q.[-\alpha_{4}]_{N}}{j_{4}}z_{2}^{j_{2}}...z_{4}^{j_{4}}$.

 \

 \

 For $h\in \{2,3,4\}$, let $r_{a}(h)$ be the largest integer $r$
 such that $B_{a}$ is divisible by $z^{r}_{h}$. We have that

 \

 \

 $r_{a}(h)= max \{0, q.\alpha_{k}+q.\alpha_{t}-(p-1)\}$\\

 where $\{k,t\}=\{2,3,4\}\setminus \{h\}$. One defines $r_{-a}(h)$
 in a similar manner.

 \

 Similarly let $r_{\pm a}(h)$ be the largest integer $\nu$ such
 that $B_{a}.B_{-a}$ is divisible by $z^{\nu}_{h}$. We have :

 \

 $r_{\pm a}(h)= q.max\{\alpha_{1}+\alpha_{h},
 \alpha_{k}+\alpha_{t}\}-(p-1)$.

 \

 \

 Now the equality $A_{a}.A_{-a}=A_{a^{\prime}}.A_{-a^{\prime}}$
 implies that $r_{\pm a}(h)=r_{\pm a^{\prime}}(h)$ and so we get
 the following equality.

 \

 \

 $\{\alpha_{1}+\alpha_{h},
 \alpha_{k}+\alpha_{t}\}=\{\alpha^{\prime}_{1}+\alpha^{\prime}_{h},
 \alpha^{\prime}_{k}+\alpha^{\prime}_{t} \}$.

 \

 \

 By an easy lemma in \cite{11} (Lemma 6.3), we conclude that there
 exists an even permutation $\sigma \in A_{4}$ of order $2$, such
 that $\alpha_{i}=\alpha^{\prime}_{\sigma(i)}$. We first claim that
 $\sigma \neq 1$. This in fact follows from the linear independence
 of the rows, see Remark 6.1, which ensures that $\alpha_{i}$ and
 $\alpha^{\prime}_{i}$ are not all equal. Furthermore, without loss
 of generality, we can assume that $\alpha_{i}=r_{1i}$ for all
 $i=1,..,4$. That is, we may consider $(\alpha_{1},...,\alpha_{4})$
 as the first row of the matrix $A$ of the abelian covering . We
 set $a_{i}=\alpha_{i}$ instead of $r_{1i}$ for simplicity. Now
 since $\alpha_{i}$ and $\alpha^{\prime}_{i}$ are different by the
 above argument, we may again, without loss of generality, suppose
 that:

 \

 $\alpha^{\prime}_{1}=a_{2}$, $\alpha^{\prime}_{2}=a_{1}$

 \

 $\alpha^{\prime}_{3}=a_{4}$, $\alpha^{\prime}_{4}=a_{3}$

 \

 \

 by our assumptions on $a_{i}$ and $a^{\prime}_{i}$, we have that

 \

 $\sum [a_{i}]_{N}=\sum [a^{\prime}_{i}]_{N}= 2N$

 \

 \

 Suppose that $[a_{1}]_{N}+[a_{2}]_{N}=[a_{3}]_{N}+[a_{4}]_{N}=N$,
 or in other words, $[a_{2}]_{N}=-[a_{1}]_{N}$ and
 $[a_{4}]_{N}=-[a_{3}]_{N}$ in $\mathbb{Z}/N\mathbb{Z}$. This means
 that the two rows $n=(a_{1},..,a_{4})$ and
 $n^{\prime}=(a^{\prime}_{1},..,a^{\prime}_{4})$ are linearly
 dependent and this contradicts the irreducibility by Remark 6.1.
 So the above equality does not hold and we may assume that
 $a_{1}+a_{2}<N$ and $a_{3}+a_{4}>N$. Now consider the row vector

 \

 $n+n^{\prime}=(a_{1},..,a_{4})+(a^{\prime}_{1},..,a^{\prime}_{4})=(a_{1}+a_{2},a_{1}+a_{2},a_{3}+a_{4},a_{3}+a_{4})$

 \

 \

 Note that Remark 6.1 assures that $n+n^{\prime} \neq \pm n$ and
 one can easily verify that this row vector also satisfies the
 conditions for $a_{i}$ and $a^{\prime}_{i}$ (in fact
 $2([a_{1}+a_{2}]+[a_{3}+a_{4}])=2([(N-1)(a_{1}+a_{2})]+[(N-1)(a_{3}+a_{4})])=2N$)
 and so we may replace the second row
 $(a^{\prime}_{1},..,a^{\prime}_{4})=(a_{2},a_{1},a_{4},a_{3})$ by
 this row vector and the equality
 $A_{n}.A_{-n}=A_{n^{\prime}}.A_{-n^{\prime}}$ should hold for this
 row vector as $n^{\prime}$ and $(a_{1},..,a_{4})$ as $n$. We show
 that this is impossible. In fact, if this equality holds, it is
 easy to see that the left hand side must contain a monomial of the
 form $z_{2}^{\alpha}z_{3}^{\beta}$ and also a monomial of the form
 $z_{1}^{\gamma}z_{4}^{\delta}$. This means that
 $a_{2}+a_{3}=a_{1}+a_{4}=a_{1}+a_{3}=a_{2}+a_{4}=N$ which is
 exactly to say that $n=n^{\prime}=(a_{1},a_{1},-a_{1},-a_{1})$.
 This is against our assumptions and this contradiction completes
 the proof. \hspace{1.5cm}$\Box$

 \

  \

 \subsection{Families with large $s$}

 \

 \

 In this section we assume that the family has Galois group of the
 form $\mathbb{Z}/n\mathbb{Z} \times \mathbb{Z}/m\mathbb{Z}$, i.e.,
 that the matrix $A$ is a $2\times s$ matrix.

 \

 \

 In \cite{11} Moonen uses the Dwork-Ogus obstruction also in order to
 prove that further examples of Shimura families of cyclic covers
 of $\mathbb{P}^{1}$ do not exist for $s> 4$. The core observation
 in his proof is then that for a minimal Shimura family not
 existing in his table, there exists two integers $n$ and
 $n^{\prime}$ such that $d_{n}=d_{-n^{\prime}}=1$ and
 $d_{-n}=d_{-n^{\prime}}=s-3$. To deduce the existence of these two
 elements, he argues that for $n\in (\mathbb{Z}/m\mathbb{Z})^{*}$,
 $d_{n}+d_{-n}=s-2$. This point does not remain necessarily true
 for families of abelian coverings. A counterexample can be given by the family\\

 $\begin{pmatrix} 1&1&1&1&0\\
 0&0&1&1&2\\
 \end{pmatrix}$ with $N=4$.\\

 For this family all of the eigenspaces are of type $(a,0)$ (or $(0,a)$) or $(1,1)$.\\

 Let us explain our strategy for the rest of this section. Recall
 from construction 2.2.2 that there exists a Shimura subvariety
 $S_f$ which contains the moduli variety $Z$ of the family $f:Y\to
 T$. As we have already explained in 2.2.2, $S_f$ is the smallest
 Shimura subvariety containing $Z$ and if $M$ is the generic
 Mumford-Tate group of the family such that $M^{ad}_{\mathbb{R}}=
 Q_{1}\times...\times Q_{l}$ is the decomposition into
 $\mathbb{R}$-simple factors, then $\dim S_f=\sum \delta(Q_{i})$
 with $\delta(Q_{i})$ as described in 2.2.2. Therefore $Z$ is a
 Shimura subvariety if and only if $\sum \delta(Q_{i})=s-3$. We are
 going to show that for large $s$ the families we consider do not
 give rise to Shimura subvarieties in $A_g$. By the above, this is
 equivalent to $\dim S_f=\sum \delta(Q_{i})>s-3$. Consider connected
 monodromy group $Mon^{0}$ of this family (see 6.2.1 below).
 Assuming on the contrary that $Z$ is a Shimura subvariety
 would imply that $Mon^0$ is a normal subgroup of $M$, see Remark 6.2.3.
 Therefore $Mon^{0,ad}_{\mathbb{R}}=\prod_{i\in K}Q_{i}$ for
 some $K\subseteq\{1,...,l\}$. In Lemma 6.2.2 we compute the
 monodromy group $Mon^{0}(\mathcal{L}_{i})$ of an eigenspace
 $\mathcal{L}_{i}$. Using this, we show that for large $s$ there
 are eigenspaces $\mathcal{L}_{j_i}$ with non-isomorphic monodromy
 groups (and hence corresponding to distinct
 $Q_i=Mon^{0,ad}_{\mathbb{R}}(\mathcal{L}_{j_i})$ for $i\in K$ in
 the above decomposition of $Mon^{0,ad}_{\mathbb{R}}$) such that $\sum
 \delta(\mathcal{L}_{j_i})>s-3$, where
 $\delta(\mathcal{L}_{j_i})=\delta(Mon^{0,ad}_{\mathbb{R}}(\mathcal{L}_{j_i}))$.
 In particular, $\dim S_f>s-3$, a contradiction. See also Remark 6.2.4.\\

 Now we state the necessary definitions and lemmas in order to
 achieve the above goals.

 \

 \

 \textbf{Definition 6.2.1.} Let $f:Y\rightarrow T$ be a family of
 abelian Galois covers of $\mathbb{P}^{1}$ as constructed in
 section $1$. Then the local system
 $\mathcal{L}=R^{1}f_{*}\mathbb{C}$ gives rise to a polarized
 variation of Hodge structures (PVHS) of weight $1$. consider the
 associated monodromy representation $\pi_{1}(T,x)\to GL(V)$, where
 $V$ is the fiber of $\mathcal{L}$ at $x$ (see for example \cite{03}, \S
 3.1.1). The Zariski closure of the image of this morphism is
 called the \emph{monodromy group} of $\mathcal{L}$. We denote the
 identity component of this group by $Mon^{0}(\mathcal{L})$. The
 PVHS decomposes according to the action of the abelian Galois
 group $G$ and the eigenspaces $\mathcal{L}_{i}$ (or
 $\mathcal{L}_{\chi}$ where $i\in G$ corresponds to character $\chi
 \in \mu_{G}$ by Remark 2.1.1) are again variations of Hodge
 structures and we are mainly interested in these. Take a $t\in T$
 and assume that $h^{1,0}((\mathcal{L}_{i})_{t})=a$ and
 $h^{0,1}((\mathcal{L}_{i})_{t})=b$. Since monodromy group respects
 the polarization of the Hodge structures ([R], 3.2.6),
 $(\mathcal{L}_{i})_{t}$ is equipped with a Hermitian form of
 signature $(a,b)$ (see \cite{85}, 2.21 and 2.23). This implies that
 $Mon^{0}(\mathcal{L}_{i}) \subseteq U(a,b)$. In this case, we say
 that $\mathcal{L}_{i}$ is \emph{of type} $(a,b)$. Two eigenspaces
 $\mathcal{L}_{i}$ and $\mathcal{L}_{j}$ of types $(a,b)$ and
 $(a^{\prime},b^{\prime})$ are said to be \emph{of
     distinct types} if $\{a,b\}\neq \{a^{\prime},b^{\prime}\}$.\\

 The above observations are key to our further analysis. Let us
 first prove a lemma which is an extension of \cite{11}, Prop 4.7 to the
 setting of abelian covers.

 \

 \

 \textbf{Lemma 6.2.2.} Let $\mathcal{L}_{i}$ be an eigenspace as
 discussed above of type $(a,b)$ with $ab\neq 0$. Then
 $Mon^{0}(\mathcal{L}_{i})=SU(a,b)$ unless when $|G|=2l$ is even
 and $i$ is of order $2$ in $G$, in which case there is a
 surjection from $Mon^{0}(\mathcal{L}_{i})$ to $SU(n,n)=Sp_{2n}$,
 where $n=d_{i}$.

 \

 \

 \textbf{Proof.} Let $t\in T$ and let $\chi\in Hom(G,
 \mathbb{C}^{*})$ be the character corresponding to $i$. Consider
 the cover $f_{\chi ,t}: Y_{\chi ,t} \rightarrow \mathbb{P}^{1}$
 with group $\chi(G)$ branched above the points $z_{j}$ with local
 monodromy $\chi(\phi(\gamma_{j}))$ about $z_{j}$, where $\phi$ is
 the surjection in Remark 2.1. Note that $\chi(G)$ is a cyclic
 group and so $f_{\chi ,t}$ is in fact a cyclic cover with group
 $\chi(G)$. Varying $t\in T$, we get a family of cyclic covers of
 $\mathbb{P}^{1}$. The eigenspace $\mathcal{L}_{i}$ is exactly the
 eigenspace corresponding to this family (or in other words, it is
 the $\mathcal{L}_{1}$ of this family of cyclic covers). Unless
 when $|G|=2l$ is even and $i$ is of order $2$ in $G$, Theorem
 5.1.1 of \cite{09} applies and we get that
 $Mon^{0}(\mathcal{L}_{i})=SU(a,b)$. If $|G|=2l$ is even and $i$ is
 of order $2$ in $G$ by taking quotient of the family $f_{\chi ,t}$, we obtain a family of hyperelliptic curves of the form $w^{2}=(z-z_{1})...(z-z_{2n+2})$. Note that in this case it follows from the formulas of Proposition 2.2.3 that there are $2n+2$ odd powers in the equation of $f_{\chi ,t}$ for $n=d_{i}$.
 Now it is well-known that $Mon^{0}$ of a versal family of
 hyperelliptic curves is the full symplectic group and so the proof
 is complete. \hspace{1.5cm}$\Box$

 \

 \

 \textbf{Remark 6.2.3.} If the family $f:Y\rightarrow T$ gives rise
 to a Shimura subvariety in $A_{g}$ then the connected monodromy
 group $Mon^{0}$ is a normal subgroup of the generic Mumford-Tate
 group $M$. In fact in this case $Mon^{0}=M^{der}$, see for example
 \cite{15} or \cite{09}. Consequently, if $M^{ad}_{\mathbb{R}}=\prod_{1}^{l} Q_{i}$ as a product of simple Lie groups then there exists a subset
 $K\subseteq\{1,...,l\}$ such that $Mon^{0,ad}_{\mathbb{R}}=\prod_{i\in K}Q_{i}$.\\

 \

 \

 \textbf{Remark 6.2.4.} Assume that $Y\rightarrow T$ is a family of
 curves and let $M$ be the generic Mumford-Tate group of this
 family. Recall from construction 2.2.2, that there is a natural
 Shimura variety $S_{f}=Sh(M,Y)$ associated to $M$ (which is a
 reductive group) and the dimension of $S_{f}$ only depends on
 $M^{ad}_{\mathbb{R}}$. The Shimura datum comes from the Hodge
 structures of the fibers in the family. This Shimura variety is
 the smallest Shimura subvariety in $A_{g}$ which contains $Z$. Our
 purpose is to show that for families of abelian covers with a
 large $s$, the moduli variety $Z$ is not a Shimura subvariety.
 This is equivalent to $\sum\delta(Q_{i})>s-3$. Here
 $\delta(Q_{i})$ is as in construction 2.2.2. Lemma 6.2.2 enables
 us to compute the connected algebraic monodromy group
 $Mon^{0}(\mathcal{L}_{i})$ of an eigenspace $\mathcal{L}_{i}$.
 Suppose that $Z$ is a Shimura subvariety and hence $\sum\delta(Q_{i})=s-3$. We show that this assumption
 leads to a contradiction. According to Remark 6.2.3 above, if $M^{ad}_{\mathbb{R}}=\prod_{1}^{l} Q_{i}$ is the decomposition
 into $\mathbb{R}$-simple factors, then $Mon^{0,ad}_{\mathbb{R}}=\prod_{i\in K}Q_{i}$ for some subset
 $K\subseteq\{1,...,l\}$. We need to find eigenspaces
 $\mathcal{L}_{j_i}$ of distinct types $\{a_i,b_i\}$ for $i\in K$
 with $a_i$ and $b_i$ large enough in the sense described below.  Concretely, if we can find eigenspaces
 $\mathcal{L}_{j_i}$ as above, then these eigenspaces, being of different types $\{a_i,b_i\}$, give rise
 to non-isomorphic $Q_i=Mon^{0,ad}_{\mathbb{R}}(\mathcal{L}_{j_i})=PSU(a_i,b_i)$ for
 $i\in K$ in the above decomposition of $Mon^{0,ad}_{\mathbb{R}}$ and if for
 $\delta(\mathcal{L}_{j_i})=\delta(Mon^{0,ad}_{\mathbb{R}}(\mathcal{L}_{j_i}))$, we have that
 $\sum \delta(\mathcal{L}_{j_i})>s-3$ (this is what we mean by \emph{large enough} in the above,
 note that $\delta(\mathcal{L}_{j_i})$ depends in our examples only on $a_i$ and $b_i$, see construction 2.2.2), then
 $\dim S_f\geq \sum \delta(\mathcal{L}_{j_i})>s-3$.

 \

 \

 \textbf{Remark 6.2.5.} An important observation is that if in the
 family one row, say the first row, does not have any $0$ entry
 then either the family of cyclic coverings arising from this row
 is a Shimura family (of cyclic covers) or the whole family (of
 abelian covers) will not be a Shimura family. This is due to the
 fact that if this family of cyclic covers is not a Shimura family
 then by the above notations and observations there are
 $\mathbb{R}$-simple factors $Q_{i}$ in the decomposition of
 $M^{ad}_{1, \mathbb{R}}$ such that $\sum\delta(Q_{i})>s-3$. Here
 $M_{1}$ is the generic Mumford-Tate group associated to this
 family of cyclic covers. As this is a sub-Hodge structure, we know
 from \cite{05} that $M_{1}$ is a quotient of $M$, the generic
 Mumford-Tate group of our family of abelian covers, and therefore
 the factors $Q_{i}$ also occur in the decomposition of $M^{ad}$
 (note that $M^{ad}$ is a semi-simple group with trivial center)
 and so $\dim S_{f}\geq \sum\delta(Q_{i})>s-3$, i.e., the family is
 not a Shimura family. On the other hand, by results of \cite{11} or
 \cite{09}, if a family of cyclic covers of $\mathbb{P}^{1}$ gives rise
 to a Shimura subvariety, then $s\leq 6$. Therefore, if $s$ is
 large, then the candidates for Shimura families are those which
 contain zeros in each row. The following theorem is meant to
 exclude these families.\\

 \

 \

 \textbf{Theorem 6.2.6.} Suppose $\sum a_{i}>2N$ and
 $\sum(-a_{i})=(\sum[-a_{i}]_{N})>2N$ and similarly $\sum b_{i}>2N$
 and $\sum(-b_{i})>2N$. For $s>19$ the family

 \

 \

 $\begin{pmatrix}
 a_{11} &\cdots& a_{1l} &a_{1l+1}&\cdots & a_{1r} &0&\cdots & 0 \\
 0  &\cdots&0&b_{1l+1}&\cdots& b_{1r}& b_{1r+1}& \cdots & b_{s} \\

 \end{pmatrix}$

 \

 \

 does not give rise to a Shimura subvariety in $A_{g}$.

 \

 \

 \textbf{Proof.} We proceed according to Remark 6.2.4, i.e., by
 showing that there are eigenspaces $\mathcal{L}_{i}$ of distinct types such that $\sum \delta(\mathcal{L}_{i})>s-3$, with $\delta(\mathcal{L}_{i})$ as in 6.2.4 and 2.2.2.\\
 Assume that the eigenspaces $\mathcal{L}_{(1,0)}$ and $\mathcal{L}_{(0,1)}$
 are respectively of types $(k_{1},r-k_{1}-2)$ and $(k_{2},s-l-k_{2}-2)$ for some
 positive integers $k_{1}$ and $k_{2}$ . If these two types are
 different (see Definition 6.2.1), we will have:

 \

 \

 $\dim S_{f}\geq k_{1}(r-k_{1}-2)+k_{2}(s-l-k_{2}-2)\geq 2(r-4)+2(s-l-4)$

 \

 \

 Now for $s\geq 19$, one sees that $2(r-4)+2(s-l-4)> s-3$
 and hence $\dim S_{f}>s-3$. It remains to treat the case where
 the two types are the same which implies that $s=r+l$. We have:

 \

 \

 $\dim S_{f}\geq k_{2}(s-l-k_{2}-2)\geq 2(s-l-4)$

 \

 \

 The right hand side is strictly greater than $s-3$ if and only if
 $r-l>5$. We may therefore assume that $r-l\leq 5$. Note that in any case, $2(s-l-4)\geq s-l-3$. We
 have:
 \

 \

 $t(-1,1)=\sum^{l}_{1}(-a_{i}) +\sum^{r}_{l+1} (b_{j}-a_{j})+
 \sum^{s}_{r} b_{k}\geq (l-1)N$.

 \

 \

 That is, for the eigenspace $\mathcal{L}_{(-1,1)}$ we have that
 $t(-1,1)\geq l-1$ and consequently $d_{(-1,1)}\geq l-2$. Similarly
 one sees that $d_{(1,-1)}\geq l-2$. Note that the eigenspaces
 $\mathcal{L}_{(-1,1)}$ and $\mathcal{L}_{(1,0)}$ (or equivalently
 $\mathcal{L}_{(0,1)}$) are of distinct types, for otherwise we
 will have $k_{2}\geq l-2$ and $s-l-k_{2}-2\geq l-2$. The first
 inequality says that $k_{2}>5$ (because $s>19$ and hence $l>7$)
 and the second one implies that $k_{2}\leq s-2l\leq 5$. This
 contradiction shows that we have an eigenspace of a new type. Now
 since $d_{(-1,1)}d_{(1,-1)}\geq (l-2)^{2}>l$, we conclude that:

 \

 \

 $\dim S_{f}\geq
 k_{2}(s-l-k_{2}-2)+d_{(-1,1)}d_{(1,-1)}>(s-l-3)+l=s-3$

 \

 and the claim follows. \hspace{1.5cm}$\Box$

 \

 \

\

\

Johannes Gutenberg-Universität Institut für Mathematik
Staudingerweg 9 55099 Mainz, Germany.

\

E-mail address: \url{mohajer@uni-mainz.de}

\

E-mail address: \url{zuok@uni-mainz.de}

\end{document}